\documentclass{amsart}

\usepackage{amssymb,amsfonts, latexsym,amsmath,hhline,array,longtable}
\usepackage{pdfsync,color, comment,colortbl, mathrsfs,stmaryrd,cite,graphicx,yfonts}

\usepackage{ifpdf}
\ifpdf \usepackage[colorlinks=true, citecolor=blue, linkcolor=blue, urlcolor=blue]{hyperref} \fi

\newcommand{\cal}{\mathcal}

\newtheorem{formula}{}[section]
\newtheorem{definition}[formula]{Definition}
\newtheorem{corollary}[formula]{Corollary}
\newtheorem{remark}[formula]{Remark}
\newtheorem{lemma}[formula]{Lemma}
\newtheorem{theorem}[formula]{Theorem}
\newtheorem*{claim}{Claim}

\def\thrm{\begin{theorem}}
\def\thrml#1{\begin{theorem}\label{#1}}
\def\ethrm{\end{theorem}}
\def\rmrk{\begin{remark}}
\def\rmrkl#1{\begin{remark}\label{#1}}
\def\ermrk{\end{remark}}
\def\dfntn{\begin{definition}}
\def\dfntnl#1{\begin{definition}\label{#1}}
\def\edfntn{\end{definition}}
\def\nmrt{\begin{enumerate}}
\def\enmrt{\end{enumerate}}
\def\tm#1{\item[{\rm (#1)}]}
\def\qtnl#1{\begin{equation}\label{#1}}
\def\eqtn{\end{equation}}
\def\lmm{\begin{lemma}}
\def\lmml#1{\begin{lemma}\label{#1}}
\def\elmm{\end{lemma}}
\def\crllr{\begin{corollary}}
\def\crllrl#1{\begin{corollary}\label{#1}}
\def\ecrllr{\end{corollary}}
\def\css{\begin{cases}}
\def\ecss{\end{cases}}
\def\prf{\begin{proof}}
\def\eprf{\end{proof}}
\def\clm{\begin{claim}}
\def\eclm{\end{claim}}

\def\cP{{\cal P}}

\def\cT{{\cal T}}

\def\cX{{\cal X}}
\def\cY{{\cal Y}}

\DeclareMathOperator{\aut}{Aut}

\DeclareMathOperator{\cyl}{Cyl}
\DeclareMathOperator{\diag}{Diag}

\DeclareMathOperator{\id}{id}

\DeclareMathOperator{\iso}{Iso}

\DeclareMathOperator{\pr}{pr}

\DeclareMathOperator{\sym}{Sym}
\DeclareMathOperator{\WL}{WL}

\def\grp#1{\langle {#1}\rangle}

\def\phmaa#1{{\phantom{x}\hspace{-2mm}^{#1}}}
\def\qaq{\quad\text{and}\quad}
\def\qoq{\quad\text{or}\quad}

\begin{document}

\title[Cartesian products of graphs and their coherent configurations]{Cartesian products of graphs and their coherent configurations}

\author{Jinzhuan Cai}
\address{Hainan University, Haikou, China}
\email{caijzh12@163.com}
\author{Jin Guo}
\address{Hainan University, Haikou, China}
\email{guojinecho@163.com}
\author{Alexander L. Gavrilyuk}
\address{University of Memphis, Tennessee, U.S.A.}
\email{a.gavrilyuk@memphis.edu}
\author{Ilia Ponomarenko}
\address{Hainan University, Haikou, China; Steklov Institute of Mathematics at St. Petersburg, Russia}
\email{inp@pdmi.ras.ru}

\thanks{}
\date{}

\begin{abstract}
The coherent configuration $\WL(X)$ of a graph $X$ is the smallest coherent configuration on the vertices of~$X$ that contains the edge set of $X$ as a relation. The aim of the paper is to study  $\WL(X)$ when $X$ is a Cartesian product of graphs. The example of a Hamming graph shows that, in general, $\WL(X)$ does not coincide with the tensor product of the coherent configurations of the factors. We prove that  if $X$ is ``closed'' with respect to the  $6$-dimensional Weisfeiler-Leman algorithm, then $\WL(X)$ is the tensor product of the coherent configurations of certain graphs related to the prime decomposition of~$X$. This  condition is trivially satisfied for almost all graphs. 
In addition, we prove that the property of a graph ``to be decomposable into a Cartesian product of $k$ connected prime graphs'' for some $k\ge 1$ is recognized by the $m$-dimensional Weisfeiler-Leman algorithm  for all $m\ge 6$.
\end{abstract}

\maketitle

\section{Introduction}
A coherent configuration can be thought of as a partition of a complete directed graph into arc-disjoint subgraphs that have “combinatorial symmetries” (precise definitions are given in Section~\ref{160224a}). The coherent configurations were introduced by B.~Weisfeiler and A.~Leman~\cite{WLe68} in connection with the graph isomorphism problem, and  independently by D.~Higman~\cite{Hig1970a} as a tool for studying permutation groups. This fact in itself hints at the intermediate position of coherent configurations between graphs and permutation groups. In fact, coherent configurations inherit the symmetries of graphs and admit a description in terms which are similar to those 
used for permutation groups~\cite{EvdP2009a}.

With each graph $X$, one can associate a uniquely determined coherent configuration $\WL(X)$ which is the output of the classical Weisfeiler-Leman algorithm~\cite{WLe68} and has the same automorphism group as~$X$. More exactly,  $\WL(X)$ is the smallest coherent configuration on the vertex set $\Omega$ of~$X$, in  which the edge set of~$X$ is a union of some  arc-disjoint subgraphs of the complete directed graph on $\Omega$. One of the main questions in this connection is which  natural properties of~$X$ can be determined from~$\WL(X)$. For example among such properties one can find connectivity, regularity, or acyclicity. On the other hand, whether the planarity of a graph~$X$ can be determined from the coherent configuration $\WL(X)$ is currently unknown, see \cite{Kiefer2017}.

In order to pose our question more precisely, we will look at it a little more generally. We consider properties of a graph that can be written as a formula in a fragment of the first-order logic with counting quantifiers, see~\cite{Grohe2017}. Those formulas that use at most $m+1$ variables can be seen inside  the output $\WL_m(X)$\footnote{It is an $m$-ary coherent configuration in the sense of L.~Babai, see~\cite{Ponomarenko2022a}.} of the $m$-dimensional Weisfeiler-Leman algorithm applied to the graph $X$, see~\cite{CaiFI1992}. The output can be thought of as a partition of the Cartesian $m$-power~$\Omega^m$, where $\Omega$ is  the vertex set  of the graph~$X$. For $m=2$, we have $\WL_m(X)=\WL(X)$, and, moreover, the natural projection $\pr_2 \WL_m(X)$ of~$\WL_m(X)$ to the Cartesian square~$\Omega^2$ is a coherent configuration which is larger than or equal to $\WL(X)$ (the definition of the partial order of coherent configuration can be found in  Subsection~\ref{300924a}). In fact,
\qtnl{300924c}
\WL(X)=\pr_2 \WL_2(X)< \pr_2 \WL_3(X)< \cdots < \pr_2 \WL_m(X)=
\cdots
\eqtn
for some positive integer $m$; here, the final term $\pr_2 \WL_m(X)$ is the coherent configuration associated with the permutation group $\aut(X)$.

Informally speaking, two graphs $X$ and $X'$ are said to be \emph{$\WL_m$-equivalent} if any formula in the fragment of first-order logic with counting quantifiers, that uses at most $m+1$ variables, is true either for both $X$ and $X'$ or for none of them (an equivalent algebraic definition of this concept is given in Subsection~\ref{260324a}).  Following \cite{Arvind2020}, a property $\cP$ of a graph is said to be \emph{$\WL_m$-invariant} if any two $\WL_m$-equivalent graphs satisfy or not the property $\cP$ simultaneously, or in other words, the $m$-dimensional Weisfeiler-Leman algorithm $\WL_m$ distinguishes any two graphs if exactly one of them satisfies~$\cP$.  In these terms, our question can be formulated as follows: for a property $\cP$ find the minimal $m$ such that $\cP$ is $\WL_m$-invariant.

In the present paper, we address the above question for the property ``to be decomposable into a nontrivial Cartesian product''.  Let $X_1=(\Omega_1,E_1)$ and $X_2=(\Omega_2,E_2)$ be graphs. The \emph{Cartesian product} $X_1\,\square\,X_2$  is a graph with vertex set $\Omega_1\times\Omega_2$, in which  two vertices $(\alpha_1,\alpha_2)$ and $(\beta_1, \beta_2)$ are adjacent precisely if 
$$
\alpha_1=\alpha_2\text{  and } (\beta_1,\beta_2)\in E_2\qoq (\alpha_1,\alpha_2)\in E_1\text{  and } \beta_1=\beta_2.
$$
The Cartesian product is one of the basic graph products and admits a unique prime factorization if the factors are connected. The proof of our first main result (Theorem~\ref{070224a} below) is based on  a well-developed theory of the Cartesian products \cite{Hammack2011}.

\thrml{070224a}
Given a positive integer $k$, the property of a graph ``to be decomposable into a Cartesian product of $k$ connected prime graphs'' is $\WL_m$-invariant for all $m\ge 6$.
\ethrm

It is not clear whether the constant $6$ in Theorem~\ref{070224a}  can be reduced. However, the following example shows that the minimal $m$ for which the statement could be true is at least~$3$. Let $X$ and $X'$ be the Hamming graph $H(2,4)$  and the Shrikhande graph (see, e.g., \cite[Example 2.6.17]{CP2019}), and let $\cP$ be the above mentioned property with $k=2$. Then $X$ satisfies $\cP$ whereas $X'$ does not. On the other hand, these graphs are strongly regular and have the same parameters; in particular, they are $\WL_2$-equivalent. Thus, the property $\cP$ is not $\WL_2$-invariant.

It would be nice to have a certificate in terms of coherent configurations for the property in Theorem~\ref{070224a}. An  example in Section~\ref{280324a} shows that if $X_1$ and~$X_2$ are complete graphs of different orders, then 
$$
\WL(X_1\,\square\,X_2)=\WL(X_1)\,\otimes\,\WL(X_2).
$$  
Thus, as such a certificate, one can consider the decomposition of the coherent configuration into a tensor product.

In general, decomposing the coherent configuration of a Cartesian product of graphs into a tensor product is not straightforward. Indeed, assume that we are given $n>1$ graphs $X_1,\ldots,X_n$ such that any two of them   are $\WL$-equivalent, and  a permutation group $G\le\sym(n)$. Then one can define a coherent configuration  that we will call the \emph{exponentiation} of the family $\{\WL(X_i)\}_{i=1}^n$
by  $G$, and denote it by $\{\WL(X_i)\}_{i=1}^n\uparrow G$  (see  Section~\ref{270324s}). If each graph $X_i$ is complete and $G=\sym(n)$, the exponentiation is equal to the coherent configuration of the Hamming graph, that is indecomposable with respect to the tensor product, though the Hamming graph  in our case is  just $X_1\,\square \cdots\square \,X_n$.  

Our second result gives a sufficient condition for the coherent configuration of a Cartesian product of graphs  to be decomposable into tensor product. To formulate the condition, we call the graph $X$ (and also its coherent configuration) \emph{$\WL_m$-closed} ($m\ge 2$)\label{220624a} if  
$$
\WL(X)=\pr_2 \WL_m(X).
$$ 
Informally, this means that the $m$-dimensional Weisfeiler-Leman algorithm applied to $X$ results in the same coherent configuration as the $2$-dimensional one. By formula~\eqref{300924c}, this is always true if the coherent configuration $\WL(X)$ is discrete (i.e., the largest coherent configuration on~$\Omega$). Since this property is satisfied for almost all graphs (see, e.g., \cite{Babai1980a}), almost every graph is $\WL_m$-closed for all $m$. Another natural example of such graphs are the Cai-F\"urer-Immerman graphs~\cite{CaiFI1992}.

\thrml{070224b}
Let $X_1,\ldots,X_n$  be  connected prime graphs. Assume that the graph $X_1\,\square \cdots\square \,X_n$ is $\WL_6$-closed, and denote by $J_1,\ldots,J_a$ the maximal subsets of  indices, such that $i,j\in J_k$ for some $1\le k\le a$ if and only if the graphs $X_i $ and $X_j$ are  $\WL$-equivalent. Then
\qtnl{120224r}
\WL(X_1\,\square \cdots\square \,X_n)=\WL(X_{J_1})\,\otimes \cdots\otimes\WL(X_{J_a})
\eqtn
where  $X_{J_k}$ is the Cartesian product of the graphs $X_j$ with $j\in J_k$,  $k=1,\ldots,a$. Moreover, for each $k$,
\qtnl{140224a}
\WL(X_{J_k})\le \{\WL(X_j)\}_{j\in J_k}\uparrow\sym(n_k) 
\eqtn
where $n_k=|J_k|$.
\ethrm

The proof of Theorem~\ref{070224b} 
in Section~\ref{280324g} 
shows that the condition that  the graph $X_1\,\square \cdots\square \,X_n$ is $\WL_6$-closed can be replaced by a weaker one. Namely, it suffices to assume that the coherent configuration of that graph is $2$-closed, see Subsection~\ref{190324g}.
 
Two natural questions arising in connection with Theorem~\ref{070224b} are as follows: (a) is it possible to replace the condition of the $\WL_6$-closedness by $\WL_m$-closedness with some $m<6$, and (b) whether the  inclusion in~\eqref{140224a} is always equality. These questions 
are also motivated by a problem 
arising from computing Delsarte's 
linear programming bound for 
sum-rank metric codes, see 
\cite{AbiadGHP2024}.
Unfortunately, in both cases, 
we have no answer; 
in particular, \cite[Problem~5.1]{AbiadGHP2024} remains unresolved.
 
To make the present paper more or less selfcontained, we recall in Section~\ref{160224a} relevant notation and  facts about coherent configurations and graphs. In Section~\ref{280324n}, we explain in some detail the concept of the $\WL_m$-equivalence and develop some tools to study $\WL_m$-equivalent graphs in an algebraic setting. The exponentiation of a family of coherent configurations by a permutation group is introduced in Section~\ref{270324s}. The coherent configurations of the Cartesian products of graphs are studied in Sections~\ref{280324a} and~\ref{270324a}. The proofs of Theorems~\ref{070224a} and~\ref{070224b} are given in Section~\ref{280324g}.

\section{Coherent configurations and graphs}\label{160224a}

\subsection{Notation.}
Throughout the paper, $\Omega$ stands for a finite set. For any $\Delta\subseteq \Omega$, we denote by~$1_\Delta$ the diagonal of the Cartesian square $\Delta\times\Delta$. The set of all classes  of an equivalence relation $e$  on a subset of $\Omega$ is  denoted by~$\Omega/e$.
 
For a binary relation $r\subseteq\Omega\times\Omega$, we set  $r^*=\{(\beta,\alpha)\colon\ (\alpha,\beta)\in r\}$. The set of all {\it neighbors} of a point $\alpha\in\Omega$ in the relation $r$ is denoted by $\alpha r=\{\beta\in\Omega\colon\ (\alpha,\beta)\in r\}$. For relations $r,s\subseteq\Omega\times\Omega$, we put 
$$
r\cdot s=\{(\alpha,\beta)\colon\ (\alpha,\gamma)\in r,\ (\gamma,\beta)\in s\text{ for some }\gamma\in\Omega\},
$$ 
which is called the \emph{dot product} of $r$ and $s$. For $\Delta,\Gamma\subseteq \Omega$, we set $r_{\Delta,\Gamma}=r\cap (\Delta\times \Gamma)$ (and abbreviate $r_{\Delta}:=r_{\Delta,\Delta}$).  For a set $S$ of relations on $\Omega$,  we denote by $S^\cup$ the set of all unions of the elements of $S$ and put $S^*=\{s^*\colon s\in S\}$ and $S^f=\{s^f\colon s\in S\}$ for any bijection $f$ from $\Omega$ to another set.
 
\subsection{Coherent configurations: basic definitions}
Let $S$ be a partition of the Cartesian square $\Omega^2$; in particular, the elements of $S$ are treated as binary relations on~$\Omega$. A pair $\cX=(\Omega,S)$ is called a {\it coherent configuration} on $\Omega$ if the following conditions are satisfied:
\nmrt
\tm{C1}  $1_{\Omega}\in S^\cup$,
\tm{C2} $S^*=S$,
\tm{C3} given $r,s,t\in S$, the number $c_{r,s}^t=|\alpha r\cap \beta s^{*}|$ does not depend on $(\alpha,\beta)\in t$. 
\enmrt

Any relation belonging to $S$ (respectively, $S^\cup$)  is called a {\it basis relation} (respectively, a {\it relation} of~$\cX$). The set of all relations is closed with respect to taking the transitive closure and the dot product. A set $\Delta \subseteq \Omega$ is called a {\it fiber} of $\cX$ if the relation $1_{\Delta}$ is basis. Any union of fibers  is called a {\it homogeneity set} of $\cX$. 

Let $\cX=(\Omega,S)$ and $\cX'=(\Omega',S')$ be two coherent configurations. A bijection $f\colon\Omega\to\Omega'$ is called an {\it isomorphism} from $\cX$ to $\cX'$ if $S^f= S'$. The isomorphism~$f$ induces a natural bijection $\varphi\colon S\to S'$, $s\mapsto s^f$. It preserves the numbers from the condition~(C3), namely, the numbers $c_{rs}^t$ and $c_{r^{\varphi},s^{\varphi}}^{t^{\varphi}}$ are equal  for all $r,s,t\in S$. Every bijection $\varphi\colon S\to S'$ having this property is called an {\it algebraic isomorphism}, written as $\varphi\colon \cX\to \cX'$. The algebraic isomorphism $\varphi\colon\cX\to\cX'$ induces a uniquely determined bijection $S^\cup\to {S'}^\cup$ denoted also by $\varphi$. We denote by $\id_\cX$ the identity algebraic isomorphism from~$\cX$ to itself.

\subsection{Coherent configurations: partial order and closure}\label{300924a}
There is a partial order\, $\le$\, of the coherent configurations on the same set~$\Omega$. Namely, given two such coherent configurations $\cX$ and $\cX'$, we set $\cX\le\cX'$ if and only if each basis relation of~$\cX$ is the union of some basis relations of~$\cX'$. The minimal and maximal elements with respect to this ordering are the {\it trivial} and {\it discrete} coherent configurations: the basis  relations of the former one are the reflexive relation $1_\Omega$ and its complement 
in $\Omega\times\Omega$ (if $|\Omega|\geq 1$), whereas the basis relations of the latter one are singletons. 

The {\it coherent closure} $\WL(r,s,\ldots)$ of some  binary relations $r,s,\ldots$ on $\Omega$, is defined to be the smallest 
coherent configuration on $\Omega$, containing each of them as a relation. When $\{r,s,\ldots\} =S(\cX)\cup\{t\}$ for some coherent configuration $\cX$ and a relation $t$,  we denote $\WL(r,s,\ldots)$ by $\WL(\cX,t)$.


\subsection{Coherent configurations: parabolics}
 An equivalence relation $e$ on a set $\Delta\subseteq\Omega$ is called a \emph{partial parabolic} of the coherent configuration~$\cX$ if $e$ is the union of some basis relations; if, in addition, $\Delta=\Omega$, then $e$ is called a \emph{parabolic} of~$\cX$. Note that the equivalence closure of any relation of $\cX$ is a partial parabolic. If $\Delta$ is a class of a partial parabolic, $e$ is a partial parabolic such that $e_\Delta\ne\varnothing$, and   $S_{\Delta/e}$ is the set of all non-empty relations $$
 s_{\Delta/e}=\{(\Lambda,\Gamma)\in (\Delta/e)^2: s_{\Delta,\Gamma}\ne\varnothing\},
 $$
 with $s\in S$, then the pair  $\cX_{\Delta/e}=(\Delta/e,S_{\Delta/e})$ is a coherent configuration called the  \emph{restriction} of~$\cX$ to~$\Delta/e$. When $e\subseteq 1_\Omega$, we identify $\Delta/e$ with $\Delta$, and put $\cX_\Delta=\cX_{\Delta/e}$ and $S_\Delta=S_{\Delta/e}$. 
   
A (partial) parabolic $e$ is said to be \emph{decomposable} if $e$ is the  union of pairwise disjoint non-empty partial parabolics; we say that $e$ is \emph{indecomposable} if it is not decomposable. Every partial parabolic is a disjoint union of uniquely determined indecomposable partial parabolics, which  are called the \emph{indecomposable components} of~$e$. 

Any algebraic isomorphism $\varphi\colon \cX\to\cX'$ induces a bijection between partial parabolics of the coherent configurations $\cX$ and $\cX'$ that preserves the property of a partial parabolic to be indecomposable \cite[Proposition~2.3.25]{CP2019}. Let  $e$ be a partial parabolic of $\cX$ and $e'=\varphi(e)$. Assume that the classes $\Delta\in\Omega/e$ and $\Delta'\in \Omega'/e'$ are such that $\varphi$ takes the indecomposable component of $e$  containing $\Delta$ as a class to the indecomposable component of $e'$  containing $\Delta'$ as a class. Then according to \cite[Exercise~2.7.31]{CP2019},  $\varphi$ induces an algebraic isomorphism 
$\varphi_{\Delta,\Delta'}\colon \cX^{}_{\Delta}\to \cX'^{}_{\Delta'}$ 
such that $\varphi_{\Delta,\Delta'}(s_{\Delta})=\varphi(s)_{\Delta'}$ 
for every $s\in S$. In particular, if $\cX=\cX'$ and $\varphi=\id_\cX$, then  so is the mapping
 \qtnl{190224a}
 \id_{\Delta,\Gamma}: S(\cX_\Delta)\to S(\cX_\Gamma),\ s_\Delta\mapsto s_\Gamma.
 \eqtn

\subsection{Coherent configurations: tensor product}
Let  $\cX_1=(\Omega_1,S_1)$ and $\cX_2=(\Omega_2,S_2)$ be  coherent configurations.  Denote by $S_1\otimes S_2$ the set of all tensor products $s_1\otimes s_2$ with $s_1\in S_1$ and  $s_2\in S_2$, where
$$
s_1\otimes s_2=\left\{((\alpha_1,\alpha_2),(\beta_1,\beta_2))\in (\Omega_1\times\Omega_2)^2:\ (\alpha_1,\beta_1)\in s_1,\ (\alpha_2,\beta_2)\in s_2\right\}.
$$
The pair $\cX_1\otimes\cX_2=(\Omega_1\times \Omega_2,S_1\otimes S_2)$ is a coherent configuration called the \emph{tensor product} of~$\cX_1$ and~$\cX_2$. The definition naturally extends to any number of factors. Note that the relations of a tensor product 
of coherent configurations are just 
the tensor products of the relations of the factors.

\subsection{Graphs}
By a graph we mean a finite simple undirected graph, i.e., a pair $X=(\Omega,E)$ of a set $\Omega$ of vertices and an irreflexive symmetric relation $E\subseteq \Omega\times\Omega$.  The elements of $E=:E(X)$ are called \emph{edges},  and $E$ is the \emph{edge set} of the graph~$X$.  Two vertices $\alpha$ and $\beta$ are adjacent (in $X$) whenever $(\alpha,\beta)\in E$ (equivalently, $(\beta,\alpha)\in E$). The subgraph of $X$ induced by a set $\Delta\subseteq \Omega$  is denoted by $X_{\Delta}=(\Delta,E_{\Delta})$. A graph $X=(\Omega,E)$ is connected if the transitive reflexive closure of $E$ equals $\Omega^2$. 
The {\it coherent configuration of a graph} $X$ is just the coherent closure of its edge set: $\WL(X)=\WL(E)$. 

\section{The Weisfeiler-Leman method}\label{280324n}

\subsection{The {\rm WL$_m$}-equivalence of graphs}\label{260324a}
Let $m\ge 1$. Given a coherent configuration  $\cX$, the $m$-dimensional Weisfeiler-Leman algorithm constructs a canonical partition $\WL_m(\cX)$ of the $m$th Cartesian power $\Omega^m$ of the underlying set~$\Omega$. The term ``canonical'' means that for any other coherent configuration $\cX'$ on $\Omega'$, every isomorphism $f\in\iso(\cX,\cX')$ naturally extends to a bijection $\Omega^m\to {\Omega'}\phmaa{m}$ taking every class of $\WL_m(\cX)$ to a class of $\WL_m(\cX')$. 

The partition $\WL_m(\cX)$ is an $m$-ary coherent configuration in the sense of L.~Babai, see~\cite{Ponomarenko2022a} and references therein. For $m=2$, the classes of the partition $\WL_2(\cX)$ coincide with $S(\cX)$ and hence ordinary coherent configurations can be treated as $2$-ary coherent configurations. In this way, the combinatorial and algebraic isomorphisms  naturally generalize to those of $m$-ary coherent configurations. In the present paper we use  $m$-ary coherent configurations and algebraic isomorphisms between them pure formally, just  to define the $\WL_m$ equivalence and to prove Theorem~\ref{190324a}; all undefined terms and facts about $m$-ary coherent configurations can be found in~\cite{Ponomarenko2022a}. 

Let  $1\le i_1<\ldots<i_k\le m$ and $K=\{i_1,\ldots,i_k\}$, where $1\le k\le m$.  Let
$$
\pr_K:\Omega^m\to\Omega^k,\ (x_1,\ldots,x_m)\mapsto(x_{i_1},\ldots,x_{i_k}),
$$ 
be the natural $K$-projection. Then $\pr_K \WL_m(\cX)=\{\pr_K(\Lambda):\ \Lambda\in\WL_m(\cX)\}$ is a $k$-ary coherent configuration and also  (see~\cite{Ponomarenko2022a})
\qtnl{240324a}
\pr_ K\WL_m(\cX)\ge \WL_k(\cX),
\eqtn 
i.e., each class of the partition on the right-hand side is a union of some classes of the  partition on the left-hand side.
Furthermore, every algebraic isomorphism $\hat\varphi:\WL_m(\cX)\to\WL_m(\cX')$ induces an algebraic isomorphism
\qtnl{240324b}
\pr_K\hat\varphi:\pr_K\WL_m(\cX)\to \pr_K\WL_m(\cX'),\ \pr_K(\Lambda)\to \pr_K(\hat\varphi(\Lambda)),
\eqtn
see ~\cite[Lemma~3.4]{Ponomarenko2022a}. When $K$ coincides with $\{1,\ldots,k\}$, we denote~$\pr_K$ and $\pr_K\varphi$ by~$\pr_k$ and~$\varphi_k$, respectively. 

Let $m\ge 2$. For any graph $X$, we put $\WL_m(X)=\WL_m(\WL(X))$. Note that $\WL_2(X)=\WL(X)$. Two graphs $X$ and $X'$ are  said to be \emph{$\WL_m$-equivalent} 
if there exists an algebraic isomorphism  $\hat\varphi:\WL_m(X)\to\WL_m(X')$ such that  $\hat\varphi_2(E)=E'$, where $E=E(X)$ and  $E'=E(X')$. One can prove that  $\WL_m$-equivalent graphs are also $\WL_k$-equivalent for $k\le m$.

 \subsection{The $2$-extensions}\label{190324g}
Let $\cX$ be a coherent configuration on $\Omega$ and $\Delta=\diag(\Omega^2)$. Following~\cite{EvdP1999c}, the \emph{$2$-extension} and  \emph{$2$-closure}  of  $\cX$  are defined to be, respectively,  the coherent configurations 
$$
\hat\cX=\WL(\cX\otimes\cX,1_\Delta) \qaq\bar \cX=(\hat\cX_\Delta)^\zeta,
$$
where $\zeta:\Delta\to\Omega$ is the mapping taking  a pair $(\alpha,\alpha)$ to the point $\alpha$, $\alpha\in \Omega$.  It holds true that $\bar\cX\ge \cX$ and if  $\bar\cX=\cX$, then  $\cX$ is said to be \emph{$2$-closed}.  Examples of $2$-closed coherent configurations include all schurian coherent configurations, i.e., those the basis relations of which are the orbits of an arbitrary permutation group acting on pairs. On the other hand, there are many coherent configurations which are not $2$-closed, see, e.g., \cite[Example~3.5.18]{CP2019}.

The basis relations of the coherent configuration $\hat \cX$ form a partition of the set~$\Omega^2\times\Omega^2$  which is naturally identified with $\Omega^4$. By \cite[Theorem~4.6.18]{CP2019}, we have
\qtnl{210324a}
\pr_4 \WL_6(\cX)\ge S(\hat\cX).
\eqtn
Furthermore, the $2$-projection of every reflexive basis relation of $\hat\cX$ is a basis relation of $\bar\cX$ (see \cite[Theorem~3.5.16]{CP2019}). Since $S(\hat\cX)$ is a partition of $\Omega^2\times\Omega^2$, this implies that $\pr_2 S(\hat\cX)\supseteq  S(\bar\cX)$. It follows that if the coherent configuration $\cX$ is $\WL_6$-closed (see page~\pageref{220624a}), then   inclusion~\eqref{210324a} implies that
$$
S(\cX)=\pr_2\WL_6(\cX) =\pr_2(\pr_4 \WL_6(\cX))\ge \pr_2 S(\hat\cX)\supseteq S(\bar\cX)\ge S(\cX).
$$
Thus, $S(\bar\cX)=S(\cX)$ and hence the coherent configuration  $\cX=\bar \cX$ is $2$-closed. Applying this statement to  $\cX=\WL(X)$, where $X$ is a $\WL_6$-closed graph,  we obtain the following lemma. 

\lmml{270324l}
Let $X$ be a $\WL_6$-closed graph. Then the coherent configuration $\WL(X)$ is $2$-closed.
\elmm

In general, it is difficult to determine the  basis relations of $\hat\cX$ explicitly. However, some of its relations are known, for example, if $s$ is a relation  of a coherent configuration~$\cX$  on $\Omega$, then
a \emph{cylindrical} relation 
of the form
\qtnl{190324h}
\cyl_s(i,j)=\{(x,y)\in \Omega^2\times \Omega^2:\ (x_i,y_j)\in s\}
\eqtn
is a relation of $\hat\cX$ for all indices $i,j\in\{1,2\}$, see \cite[Theorem~3.5.7]{CP2019}. The proof of the following theorem is based on some facts about $m$-ary coherent configurations; all of them can be found in~\cite{Ponomarenko2022a} and~\cite{Chen2023}.

\thrml{190324a}
Let $X$ and $X'$ be two $\WL_m$-equivalent graphs, $m\ge 6$ an integer, and let $\cX=\WL(X)$ and $\cX'=\WL(X')$. Then there are algebraic isomorphisms $\varphi:\cX\to\cX'$ and  $\hat\varphi:\hat\cX\to\hat\cX'$ such that 
\qtnl{190324j}
\varphi(E(X))=E(X')\qaq \hat\varphi(\cyl_s(i,j))=\cyl_{\varphi(s)}(i,j)
\eqtn
for all $s\in S(\cX)^\cup$ and all $i,j\in\{1,2\}$.
\ethrm
\prf
By the assumption of the theorem, there exists an algebraic isomorphism  $\hat\psi:\WL_m(\cX)\to\WL_m(\cX')$ such that  $\hat\psi_2(E)=E'$, where $E=E(X)$ and  $E'=E(X')$. Put $\varphi=\hat\psi_2$. Then the first equality in~\eqref{190324j} is trivial. 

To prove the second statement,  we assume, without loss of generality,  that $m=6$. Then  from \cite[Lemma~4.6.17]{CP2019}, it follows that $\pr_4 \WL_m(\cX)=S(\tilde\cX)$  for a certain coherent configuration $\tilde\cX$ on $\Omega^2$. Similarly, $\pr_4 \WL_m(\cX')=S(\tilde\cX')$  for a certain coherent configuration $\tilde\cX'$ on ${\Omega'}^2$. We claim that the bijection
$$
\hat\psi_4:S(\tilde\cX)\to S(\tilde\cX')
$$ 
is an algebraic isomorphism from $\tilde\cX$ to $\tilde\cX'$. Indeed, from the proofs of \cite[Lem\-ma~4.6.17]{CP2019} and  \cite[Lem\-ma~3.5]{Chen2023}, it follows that the intersection numbers of the coherent configurations $\tilde\cX$ and $\tilde\cX'$ can be expressed via the intersection numbers of the $m$-ary coherent configurations $\WL_m(\cX)$ and $\WL_m(\cX')$. Thus the claim follows from the fact that $\hat\psi$ is an algebraic isomorphism.

The algebraic isomorphism $\pr_4\hat\psi$ preserves the cylindrical relations. Namely, let $s$ be a relation of the coherent configuration $\cX$ and $i,j\in\{1,2\}$. Then
\qtnl{240324j}
\cyl_s(i,j)\in S(\tilde\cX)^\cup\qaq \hat\psi_4(\cyl_s(i,j))=\cyl_{\varphi(s)}(i,j).
\eqtn
Indeed, $\cyl_s(i,j)=\pr_K^{-1}(s)$, where $K=\{i,j+2\}$ and  $\pr_K$ is the projection from~$\Omega^4$ to $\Omega^2$. This proves the first equality in~\eqref{240324j}. Next, from formula~\eqref{240324b}, it follows that a class $\Lambda\in\WL_4(\cX)$ is contained in  $\pr_K^{-1}(s)$ if and only if the class $\hat\psi_4(\Lambda)$ is contained in 
$$
\pr_K^{-1}(\hat\psi_2(s))=\pr_K^{-1}(\varphi(s))=\cyl_{\varphi(s)}(i,j).
$$ 
This proves the second equality in formula~\eqref{240324j}.

Note that $1_\Delta$ is equal to the  intersection of the relations $\cyl_s(i,j)$ with $s=1_\Omega$,  and   also $
s_1\otimes s_2=\cyl_{s_1}(1,1)\cap \cyl_{s_2}(2,2)$  for all  basis relations $s_1$ and $s_2$ of $\cX$. By the definition of~$\hat\cX$ and formula~\eqref{240324j}, this implies that 
$$
\tilde\cX\ge \WL(\cX\otimes\cX,1_\Delta)=\hat\cX\qaq \hat\psi_4(\tilde\cX)=\tilde\cX'.
$$ 
It follows that  the restriction $\hat\varphi$ of the algebraic isomorphism $\hat\psi$ to the coherent configuration $\hat\cX$ is an algebraic isomorphism from $\hat\cX$ to $\hat\cX'$. Then the second equality in~\eqref{190324j} follows from~\eqref{240324j}.
\eprf
 
\section{A generalization of exponentiation}\label{270324s}
 
Let $n\ge 1$ be an  integer, $\cX_i$ a coherent configuration on~$\Omega_i$,  and  $\varphi_{ij}:\cX_i\to \cX_j$  an algebraic isomorphism, $1\le i,j\le n$. We assume that for all $i,j,k$, 
\qtnl{170324a}
\varphi_{ii}=\id_{\cX_i}\qaq \varphi_{ij}\varphi_{jk}=\varphi_{ik}.
\eqtn
In particular, $\varphi^{}_{ij}=\varphi_{ji}^{-1}$ for all $i,j$. Let $G\le\sym(n)$ be a permutation group. For each basis relation $s_1\otimes\ldots\otimes s_n$ of the  coherent configuration $\cX=\cX_1\otimes\ldots\otimes\cX_n$ and for each permutation $g\in G$, denote by~$\varphi_g$ a permutation of the basis relations of~$\cX$, such that 
\qtnl{110324g}
(s_1\otimes\cdots\otimes s_n) ^{\varphi_g}=\varphi_{j_11}(s_{j_1})\otimes\cdots\otimes \varphi_{j_nn}(s_{j_n}),
\eqtn 
where $j_1=1^{g^{-1}},\ldots, j_n=n^{g^{-1}}$.

\lmml{130324a} 
For each $g\in G$, the bijection $\varphi_g$ is an algebraic automorphism of the coherent configuration $\cX$. 
\elmm
\prf
Let $r=r_1\otimes\ldots\otimes r_n$, $s=s_1\otimes\ldots\otimes s_n$, and $t=t_1\otimes\ldots\otimes t_n$ be basis relations of the coherent configuration~$\cX$.  Put $r'=\varphi_g(r)=r'_1\otimes\ldots\otimes r'_n$, $s'=\varphi_g(s)=s'_1\otimes\ldots\otimes s'_n$, and $t'=\varphi_g(t)=t'_1\otimes\ldots\otimes t'_n$. According to \cite[formula (3.2.7)]{CP2019}, we have 
$$
c_{r',s'}^{t'}=\prod_{i=1}^n c_{r'_i,s'_i}^{t'_i}.
$$
By formula~\eqref{110324g},  each factor of the above product can be written as 
$$
c_{r'_i,s'_i}^{t'_i}=c_{\varphi_{j_ii}(r_{j_i}),\varphi_{j_ii}(s_{j_i})}^{\varphi_{j_ii}(t_{j_i})}=c_{r_{j_i},s_{j_i}}^{t_{j_i}}.
$$
Using the same formulas for $r$, $s$ and $t$, we obtain
$$
c_{r',s'}^{t'}=\prod_{i=1}^n c_{r'_i,s'_i}^{t'_i}=\prod_{i=1}^n c_{r_{j_i},s_{j_i}}^{t_{j_i}}=
\prod_{i=1}^n c_{r_i,s_i}^{t_i}=c_{r,s}^t,
$$
as required.
\eprf

Now we define the \emph{exponentiation} $\{\cX_i\}_{i=1}^n\uparrow G$ of the family of the coherent configurations $\cX_i$, $1\le i\le n$, by a permutation group $G\le\sym(n)$ to be  the algebraic fusion of the tensor product $\cX_1\otimes\ldots\otimes\cX_n$ with respect to the group $\{\varphi_g:\ g\in G\}$. Thus the exponentiation is a coherent configuration each basis relation of  which is of the form 
\qtnl{180324z}
(s_1\otimes \cdots\otimes s_n)^G=\bigcup_{g\in G}(s_1\otimes \cdots\otimes s_n)^{\varphi_g}
\eqtn
for some basis relations $s_1,\ldots,s_n$ of the coherent configurations $\cX_1,\ldots,\cX_n$, respectively.  In the special case  $\cX_i=\cX_j$  and $\varphi_{ij}=\id_{\cX_i}$ for all $i,j$, our construction is just the exponentiation of $\cX_1$ by $G$, see \cite[Subsection~3.4.2]{CP2019}.

\section{Cartesian product of graphs}\label{280324a}

Let $I=\{1,\ldots,n\}$, and let  $X_i=(\Omega_i,E_i)$ be a graph, $i\in I$. Following~\cite{Hammack2011}, the \emph{Caresian product} 
\qtnl{080224a}
X=X_1\square\, \cdots\square\, X_n
\eqtn 
of $X_1,\ldots,X_n$ is defined to be the graph with vertex set $\Omega=\Omega_1\times\ldots\times\Omega_n$ and edge set $E=c_1\cup\cdots\cup c_n$, where for each $i\in I$ we set
\qtnl{080224b}
c_i=c_i(X)=1_{\Omega_1}\otimes\ldots\otimes 1_{\Omega_{i-1}}\otimes E_i\otimes
1_{\Omega_{i+1}}\otimes\ldots\otimes 1_{\Omega_n}.
\eqtn
In what follows, we always assume that the graph $X_i$ is connected for each~$i$. In this case, in a sense of~\cite{Chen2021},  the equivalence relations $\grp{c_i}$, $i\in I$, form the standard atomic Cartesian decomposition of a set $\Omega$.

The graph $X$ is called \emph{prime} if whenever equality \eqref{080224a} holds for some graphs $X_1, \ldots,X_n$, each with at least two vertices, it follows that $n=1$. If the graphs $X_1, \ldots,X_n$ are prime, then  we refer to equality \eqref{080224a} as a \emph{prime decomposition} of~$X$. Every connected graph admits a unique prime decomposition  up to isomorphisms and order of the factors \cite[Theorem~6.6]{Hammack2011}.

Clearly,  $c_i$ is a relation of the coherent configuration $\WL(X_1)\otimes\ldots\otimes \WL(X_n)$ for each~$i\in I$, Therefore so is the relation~$E=c_1\cup\cdots\cup c_n$. By the minimality of the coherent closure $\WL(X)=\WL(E)$, we conclude that
\qtnl{031024a}
\WL(X)\le\WL(X_1)\otimes\ldots\otimes \WL(X_n).
\eqtn
In the following statement, we give a necessary and sufficient condition for equality in this formula.

\lmml{080224d}
Let $X$ be the Cartesian product~\eqref{080224a} of connected graphs. Then
\qtnl{150324a}
\WL(X)=\WL(X_1)\otimes\ldots\otimes \WL(X_n)
\eqtn
if and only if the equivalence relation $\grp{c_i}$ is a parabolic of the coherent configuration $\WL(X)$ for each $i\in I$. 
\elmm
\prf
Assume that equality \eqref{150324a} holds. Let $i\in I$. Since the graph $X_i$  is connected, $\grp{E_i}=\Omega_i^2$. By  formula~\eqref {080224b}, this implies that 
$$
\grp{c_i}=1_{\Omega_1}\otimes\cdots\otimes 1_{\Omega_{i-1}}\otimes \Omega_i^2\otimes
1_{\Omega_{i+1}}\otimes\cdots\otimes 1_{\Omega_n}
$$
Since the right-hand side is obviously a relation of $\cX=\WL(X)$, the equivalence relation~$\grp{c_i}$ is a parabolic of $\cX$.

Conversely, assume that $e_i=\grp{c_i}$ is a parabolic of the coherent configuration~$\cX$ for each $i\in I$. In view of formula~\eqref{031024a}, it suffices to verify that 
\qtnl{160324r}
\cX\ge \WL(X_1)\otimes\cdots\otimes \WL(X_n).
\eqtn
First, we note that $E\cap e_i=c_i$ is a relation of $\cX$ for each~$i\in I$.  It is not hard to verify that then $\cX$ is larger than or equal to
$$
\WL(c_i)=\WL(\ldots\otimes  1_{\Omega_{i-1}}\otimes E_i\otimes
1_{\Omega_{i+1}}\otimes\ldots)\ge
\ldots\otimes  \cT_{\Omega_{i-1}}\otimes \WL(E_i)\otimes
\cT_{\Omega_{i+1}}\otimes\ldots.
$$
This implies that given a basis relation $s_i$ of $\cX_i$, $i\in I$,  the coherent configuration~$\cX$ contains the relation $\bar s_i=\ldots\otimes 1_{\Omega_{i-1}}\otimes s_i\otimes 1_{\Omega_{i+1}}\otimes \ldots$. Therefore, 
$$
s_1\otimes s_2\otimes\cdots\otimes s_n=\bar s_1\cdot \bar s_2\cdot \ldots \cdot \bar s_n.
$$
is also a relation of~$\cX$. This proves formula~\eqref{160324r}.
\eprf

{\bf Example.}\footnote{A more general statement was proved in~\cite{AbiadGHP2024}.}\label{270324c} Let $X_1$ and $X_2$ be complete graphs of  different orders $n_1\ge 3$ and $n_2\ge 3$, respectively. Then $\WL(X_1\,\square\,X_2)=\WL(X_1)\,\otimes\,\WL(X_2)$. To prove this equality, denote by  $A_i$ the adjacency matrix of  the relation  $c_i$, $i=1,2$. Then 
$$
(A_i)^2=(n_i-2)A_i+(n_i-1)A_0\qaq A_iA_j=A_jA_i,\qquad i,j=1,2,
$$ 
where $A_0$ is the identity matrix of order $n_1n_2$. Furthermore, the adjacency matrix of the graph $X$ is equal to $A=A_1+A_2$ and
$$
A^2=A_1^2+2A_1A_2+A_2^2=(n_1-2)A_1+(n_2-2)A_2+A',
$$
where $A'$ is a matrix such that  $A'\circ A=0$ (here, $\circ$ denotes the Hadamard multiplication). According to the Wielandt principle (see \cite[Theorem~2.3.10]{CP2019}), the set $s_i$ of all pairs $(\alpha,\beta)\in\Omega^2$ such that $(A^2\circ A)_{\alpha,\beta}=n_i-2$ is a relation of the coherent configuration $\WL(X_1\,\square\,X_2)$ for each~$i$. On the other hand, $s_i=c_i$, because $n_1\ne n_2$ by the assumption. Thus, $\grp{c_i}$ is a parabolic of $\WL(X_1\,\square\,X_2)$  and we are done by Lemma~\ref{080224d}.\hfill$\square$\medskip

When the graphs in the Cartesian product~\eqref{080224a}  are pairwise $\WL$-equivalent, the coherent configuration $\WL(X)$ cannot be a tensor product. An easy example is given by the Cartesian product of complete graphs of the same order~$k$. In this case,  $X$ is a Hamming graph and the coherent configuration $\WL(X)$ is a Hamming scheme. If $k\ge 3$, then the latter scheme has only trivial parabolics and hence cannot be a nontrivial tensor product.

\lmml{080224e}
Let $X$ be the Cartesian product~\eqref{080224a} of connected graphs. Assume that the graphs $X_1, \ldots, X_n$ are pairwise $\WL$-equivalent. Then
\qtnl{160324v}
\{\WL(X_1),\ldots,\WL(X_n)\}\uparrow{\sym(n)}\ge \WL(X).
\eqtn
\elmm
\prf
 Denote by $\cY$ the coherent configuration on the left-hand side of the inclusion~\eqref{160324v}. It suffices to verify that $E=E(X)$ is a relation of~$\cY$. Let $I=\{1,\ldots,n\}$. By  assumption,  for any $i,j\in I$, there is an algebraic isomorphism $\varphi_{ij}:\WL(X_i)\to \WL(X_j)$ such that $\varphi_{ij}(E_i)=E_j$, and also conditions~ \eqref{170324a} are satisfied for all $k\in I$. Next, $E_1$ being a relation of $\WL(X_1)$ is a union of its basis relations, say $s_{11},\ldots,s_{1a}$ for some $a\ge 1$. Then $E_i$ is a union of basis relations $s_{ij}:=\varphi_{1i}(s_{1j})$, $j=1,\ldots,a$. Since $X$ is the Cartesian product of $X_1,\ldots,X_n$, we have
 $$
 E=\bigcup _{i=1}^nc_i=\bigcup_{i=1}^n\bigcup_{j=1}^a \hat s_{ij}= \bigcup_{j=1}^a\left[\bigcup_{i=1}^n \hat s_{ij}\right],
 $$
 where $\hat s_{ij}= 1_{\Omega_1}\otimes\cdots\otimes 1_{\Omega_{j-1}}\otimes s_{ij}\otimes 1_{\Omega_{j+1}}\otimes\cdots\otimes 1_{\Omega_n}$ is a relation of the tensor product $\WL(X_1)\otimes\cdots\otimes\WL(X_n)$. It remains to verify that $\hat s_{1j}\cup\cdots \cup \hat s_{nj}$ is a relation of $\cY$ for each $1\le j\le  a$. To this end,  set $G=\sym(n)$ and denote by~$G_i$, where $i\in I$, the coset of all permutations $g\in G$ such that $1^g=i$. Then  by formula~\eqref{110324g}, we have  $(\hat s_{1j})^{\varphi_g}=\hat s_{ij}$ for all $g\in G_i$. It follows by~\eqref{180324z} that
 $$
(\hat s_{1j})^G=\bigcup_{g\in G}(\hat s_{1j})^{\varphi_g}= \bigcup_{i=1}^n\bigcup_{g\in G_i}(\hat s_{1j})^{\varphi_g}=\hat s_{1j}\cup\cdots \cup \hat s_{nj}.
$$
This completes the proof, because $(\hat s_{1j})^G$ is a relation of~$\cY$.
\eprf

\section{The $2$-extension of Cartesian product}\label{270324a}
Let  $X_i$, $i\in I$, be  connected graphs. The \emph{product relation} $c(X)$ associated with their Cartesian product~\eqref{080224a} is a relation on $E$ that consists of all pairs $(u,v)\in E\times E$ for which there exists $i\in I$ such that $u,v\in c_i(X)$, or, equivalently, 
$$
c(X)=c_1(X)^2\,\cup\ldots\cup\,c_n(X)^2.
$$
Thus, the number of classes of the equivalence relation $c(X)$ is equal to the number of factors of the Cartesian product~\eqref{080224a}.

The product relation admits an explicit description in the case when decomposition~\eqref{080224a} is prime. Namely, in accordance with \cite[Theorem~23.2]{Hammack2011}, we have 
\qtnl{190324e}
c(X)=\grp{\Theta(X)\cup\tau(X)}
\eqtn
where  $\Theta=\Theta(X)$ and $\tau=\tau(X)$ are relations  on the edge set~$E$ of the graph~$X$ that are defined  as follows (see \cite[Section~23.1]{Hammack2011}): two pairs $(x,y),(x',y')\in E$ belong to the relation
\nmrt
\tm{1}  $\Theta$ if and only if $\partial(x,y') +\partial(y,x') \ne \partial(x,x') + \partial(y,y')$,
\tm{2}  $\tau$ if and only if  $\{x\}=\{x'\}=yE\cap y'E$,
\enmrt
where $\partial(\cdot,\cdot)$ is the distance function of the graph $X$. In turn, these relations can easily be described in terms of the relations~\eqref{190324h} and some relations of the coherent configuration  $\cX=\WL(X)$. Indeed,  denote by $s_i=s_i(X)$ the set of all pairs $(\alpha,\beta)$ at distance $i$ in the graph~$X$. Then
\qtnl{190324w}
\Theta=\bigcap_{a+b\ne c+d} s(a,b,c,d),
\eqtn
where given nonnegative integers $a,b,c,d$, we put
$$
s(a,b,c,d)=s(a,b,c,d;X)=\cyl_{s_a}(1,2)\cap\cyl_{s_b}(2,1)\cap\cyl_{s_c}(1,1)\cap\cyl_{s_d}(2,2).
$$
Furthermore, set $s'=s'(X)$ to be the union of all relations $t\in S(\cX)$ for which the sum $c_{r,s}^t$ over all $r,s\in S(\cX)$ such that $r,s\subseteq E$, is equal to $1$. Then
\qtnl{190324q}
\tau=\cyl_{1_\Omega}(1,1)\cap \cyl_{s'}(2,2).
\eqtn

\lmml{090224a}
If the decomposition \eqref{080224a} is prime, then  $\Theta(X)$ and $\tau(X)$ are relations of the coherent configuration $\hat\cX$ and, in particular,  $c(X)$ is a partial parabolic of $\hat\cX$.
\elmm
\prf
The last statement follows from the first one by formula~\eqref{190324e}. Recall that $s_i\in S(\cX)^\cup$ for all nonnegative integers~$d$, see \cite[Theorem~2.6.7]{CP2019}. It follows that $s(a,b,c,d)$ is a relation of $\hat\cX$ for all $a,b,c,d$, see Section~\ref{190324g}. Thus, $\Theta$ is also a relation of~$\hat\cX$ by formula~\eqref{190324w}. Similarly, $s'$ is obviously a  relation of~$\cX$, and hence~$\tau$ is  a relation of~$\hat\cX$ by formula~\eqref{190324q}. 
\eprf

In what follows, we always assume that the decomposition \eqref{080224a} is prime.  Denote by $e$ the equivalence relation $\grp{s}\subseteq E^2$, where $s=(e_1\cap c)\cup (e_2\cap c)$ with $e_1=\cyl_{1_\Omega}(1,1)$,  $e_2=\cyl_{1_\Omega}(2,2)$, and $c=c(X)$. Clearly, $e$ is a partial parabolic of the coherent configuration~$\hat\cX$. It is not hard to see that each class of the partial parabolic~$e$ is of the form 
\qtnl{120224a}
\Delta(\bar x,i)=1_{x_1}\otimes \cdots\otimes 1_{x_{i-1}}\otimes E_i\otimes  1_{x_{i+1}}\otimes\cdots\otimes 1_{x_n},
\eqtn
for some $(n-1)$-tuple $\bar x=(x_1,\ldots,x_{i-1},x_{i+1},\ldots,x_n)$ and $i\in I$. 

\lmml{120224b}
Let $\Delta=\Delta(\bar x,i)$ and $\Gamma=\Delta(\bar y,j)$ be classes of the partial parabolic~$e$. Assume that  the graphs $X_i$ and $X_j$ are not $\WL$-equivalent. Then $\Delta$ and $\Gamma$ are classes of different indecomposable  components of~$e$.
\elmm
\prf
Assume on the contrary that $\Delta$ and $\Gamma$ are classes of some indecomposable  component of~$e$. Let $\hat\varphi:\hat\cX_\Delta\to\hat\cX_\Gamma$ be the algebraic isomorphism  defined by formula~\eqref{190224a}. Then $\hat\varphi$ induces the algebraic isomorphism 
$\psi:=\hat\varphi_{\Delta/e_1,\Gamma/e_1}$ from $\hat\cX_{\Delta/e_1}$ to~$\hat\cX_{\Gamma/e_1}$, where (the parabolic) $e_1=\cyl_{1_\Omega}(1,1)$ is as above.

On the other hand,  the classes of the equivalence relation $(e_1)_\Delta$ have the form
$$
\bar\alpha =1_{x_1}\otimes \cdots\otimes 1_{x_{i-1}}\otimes (\{\alpha\}\times \alpha E_i)\otimes  1_{x_{i+1}}\otimes\cdots\otimes 1_{x_n},
$$
where $\alpha$ runs over $\Omega_i$. It follows that there is a bijection $f:\Omega_i\to {\Delta/e_1}$ that takes $\alpha$ to $\bar\alpha$. This bijection takes $E_i$ to the set of all pairs $(\bar\alpha,\bar\beta)\in (\Delta/e_1)^2$ such that $(\alpha,\beta)\in E_i$. It is straightforward to check that this set coincides with the restriction of the relation $E\otimes\Omega^2\in S(\hat\cX)^\cup$ to the set $\Delta/e_1$, i.e.,
$$
(E_i)^f=(E\otimes\Omega^2)_{\Delta/e_1}.
$$
The algebraic isomorphism $\psi$ takes the relation on the right-hand side to the relation $(E\otimes\Omega^2)_{\Gamma/e_1}$. Repeating the above argument for the class of $(e_1)_\Gamma$, we find a bijection  $g:\Omega_j\to {\Gamma/e_1}$ such that $(E_j)^g=(E\otimes\Omega^2)_{\Gamma/e_1}$. Now, let 
$$
\varphi_f:\WL(E_i)\to \WL(E_i^f),\ s\to s^f, \qaq 
\varphi_g:\WL(E_j)\to \WL(E_j^g), \ s\to s^g,
$$
be the algebraic isomorphisms induced by the bijections $f$ and  $g$.
Then the composition $\varphi_f\,\psi\,(\varphi_g)^{-1}$ is also an algebraic isomophism from the coherent configuration $\WL(E_i)=\WL(X_i)$ to the coherent configuration $\WL(E_j)=\WL(X_j)$ that takes $E_i$ to $E_j$. Consequently, the graphs $X_i$ and~$X_j$ are $\WL$-equivalent, a contradiction.
\eprf

\section{Proof of Theorems~\ref{070224a} and ~\ref{070224b}}\label{280324g}

\subsection{Proof of Theorem~\ref{070224a}} Let $X$ and $X'$ be  $\WL_m$-equivalent graphs, $m\ge 6$.  In what follows, $\cX=\WL(X)$ and $\cX'=\WL(X')$, and $\varphi$ and $\hat\varphi$ are the algebraic isomorphisms from Lemma~\ref{190324a}. We need to verify that if one of the graphs, say~$X$, is a Cartesian product of $k$ connected prime graphs, then so is the other one, or equivalently, that the product relations $c(X)$ and $c(X')$ have the same number of classes.

By the second part of Lemma~\ref{090224a}, the product relation $c(X)$ is a partial parabolic of the coherent configuration $\hat\cX$. Since any partial parabolic has  the same number of classes as its image under an algebraic isomorphism, it suffices to verify that  
\qtnl{190324f}
\hat\varphi(c(X))=c(X').
\eqtn
By the first part of Lemma~\ref{090224a} and formula~\eqref{190324e}, this is true if $\hat\varphi(\Theta(X))=\Theta(X')$ and $\hat\varphi(\tau(X))=\tau(X')$. On the other hand,  the first equality in~\eqref{190324j} implies  that $\varphi(s_i(X))=s_i(X')$ for all nonnegative integers~$i$, see~\cite[Exercise~2.7.55(1)]{CP2019}. By the second equality in~\eqref{190324j}, this yields $\hat\varphi(s(a,b,c,d;X))=s(a,b,c,d;X')$ for all nonnegative integers~$a,b,c,d$. Thus by formula~\eqref{190324w}, we obtain
$$
\hat\varphi(\Theta(X))=\hat\varphi(\bigcap_{a+b\ne c+d} s(a,b,c,d;X))=\bigcap_{a+b\ne c+d} s(a,b,c,d;X')=\Theta(X').
$$
Similarly, $\varphi(s'(X))=s'(X')$ and using again the second equality in~\eqref{190324j}, and formula~\eqref{190324q}, we conclude that 
$$
\hat\varphi(\tau(X))=\hat\varphi(\cyl_{1_\Omega}(1,1)\cap \cyl_{s'(X)}(2,2))=
\cyl_{1_\Omega}(1,1)\cap \cyl_{s'(X')}(2,2)=\tau(X'),
$$
as required.

\subsection{Proof of Theorem~\ref{070224b}}
We keep the notation of Section~\ref{270324a}. To prove equality~\eqref{120224r}, put $Y_k=X_{J_k}$, $1\le k\le a$.  Then $X=Y_1\,\square \ldots\square \,Y_a$.  Denote by $\hat e_k$ a partial equivalence relation on $\Omega^2$ with classes $\Delta(\bar x,j)$, where $\bar x$ is as in formula~\eqref{120224a} and $j\in J_k$. Then obviously
$$
e=\hat e_1\,\cup\cdots\cup\,\hat e_a.
$$
Let $1\le k\le a$. By the definition of the set $J_k$, for no two indices $j\in J_k$ and $j'\not\in J_k$, the graphs $X_{j^{}}$ and~$X_{j'}$ are $\WL$-equivalent. By Lemma~\ref{120224b}, this implies  that $\hat e_k$ is a union of indecomposable components of the parabolic~$e$. Consequently, the support $\Omega(\hat e_k)$ of $\hat e_k$ is a homogeneity set of the coherent configuration~$\hat\cX$. By \cite[Theorem~3.5.16]{CP2019}, this implies that $\Omega(\hat e_k)$ is a relation of the coherent configuration~$\bar\cX$, and hence of the coherent configuration $\cX=\bar\cX$ (recall that $\cX$ is $2$-closed by Lemma~\ref{270324l}).  On the other hand,
$$
\Omega(\hat e_k)=\bigcup_{j\in J_k} c_j(X):=c_{J_k}(X).
$$
Thus, $e'_k:=\grp{c_{J_k}(X)}$ is a parabolic of the coherent configuration $\cX$. It is not hard to see that the parabolic $e'_k$ coincides with the  parabolic $\grp{c_k(X)}$ with respect to the Cartesian product $X=Y_1\,\square \cdots\square\, Y_a$. Since this is true for all $k$, we conclude by Lemma~\ref{080224d}  that
$$
\WL(X)=\WL(Y_1)\,\otimes \ldots\otimes \,\WL(Y_a).
$$
as required. The inclusion in \eqref{140224a} immediately follows from Lemma~\ref{080224e}.

\bibliographystyle{amsplain}

\end{document}